\documentclass[11pt]{article}
\usepackage{amssymb}
\font\tenBbb=msbm10
\def\Z{\hbox{\tenBbb Z}}
\def\T{{\hbox{\tenBbb T}}}
\font\tenBbb=msbm10

\def\sgn{\hbox{sgn}}
\font\tenBbb=msbm10

\font\tenBbb=msbm10
\newtheorem{theorem}{Theorem}[section]
\newtheorem{pro}{Proposition}[section]
\newtheorem{lem}{Lemma}[section]

\newtheorem{rem}{Remark}[section]
\newcommand{\R}{{\Bbb R}}

\def\supp{\mathop{\rm supp}\nolimits}

\def\into{\int \hspace*{-4mm} - \,}

\newcommand{\re}[1]{(\ref{#1})}

\begin{document}
\begin{center}
\noindent {\Large \bf{Well-posedness in $ H^1 $ for the  (generalized)  Benjamin-Ono equation
 on the circle}}
\end{center}
\vskip0.2cm
\begin{center}
\noindent
{\bf Luc Molinet}\\
{\small L.A.G.A., Institut Galil\'ee, Universit\'e Paris-Nord,\\
93430 Villetaneuse, France.} \vskip0.3cm
{\bf Francis Ribaud}\\
{\small Universit\'e de Marne--La--Vall\'ee, Equipe d'Analyse et de
Math\'ematiques Appliqu\'ees,\\
5 bd. Descartes, Cit\'e Descartes, Champs-sur-Marne, \\ 77454 Marne-La-Vall\'ee Cedex 2, France.} \vskip0.3cm \noindent
E-mail : molinet@math.u-paris13.fr , ribaud@math.univ-mlv.fr
\end{center}
\vskip0.5cm \noindent {\bf Abstract.} {\small We prove the local well posedness of  the Benjamin-Ono equation
and the generalized Benjamin-Ono equation in $ H^1(\T) $. This leads to a global well-posedness result in
 $ H^1(\T)$  for the Benjamin-Ono equation.}

\
\section{Introduction, main results and notations}
\subsection{Introduction and main results}
In this paper we study  the $H^1(\T)$ local well-posedness problem for  the Benjamin-Ono equation and the generalized Benjamin-Ono equation
$$
 \left\{\begin{array}{llll}
\partial_tu +{\cal{H}}\partial_{xx} u = u^k \partial_xu \,  \, ,\; (t,x)\in \R\times \T\;,\\
u(0,x)=u_0(x) \; .
\end{array}\right. \leqno{\mbox{(GBO)}}
$$
Here $\T=\R/(2\pi \Z)$, $k \geq 1$ is an integer and ${\cal{H}}$
denotes the Hilbert transform defined for  $ 2\pi \lambda $-periodic
functions by
$$
 \left\{\begin{array}{llll}
{\widehat {\cal{H}}({f})}(q)= -i \,\sgn(q)
\hat{f}(q) \, , \; q \in {\lambda^{-1}\Z^*} \, , \\
{\widehat {\cal{H}}({f})}(0)=0  \; .
\end{array}\right.
$$
\vskip0.3cm

When $k=1$, $(GBO)$ is the well known Benjamin-Ono equation. This equation  has been derived as a model for the propagation of long internal gravity waves in deep and stratified fluids \cite{B}, and, at least when the spatial domain is the whole real line, has been studied in a large amount of works in the last decades. The Benjamin-Ono equation is a totally integrable system \cite{CW}, and possesses among others the three following invariant
quantities :
$$
 I(u)=\int \, u(t,x) \, dx \,,\quad M(u)=\int \, u^2(t,x) \, dx \,,
$$
and
$$
  F(u)=\int \, u_x^2(t,x) -\frac{3}{4} \, u^2(t,x) {\cal{H}}u_x(t,x) -\frac{1}{8} \, u^4(t,x) \, dx  \;\;.
$$
Using Sobolev embedding theorems together with standard
interpolation inequalities,  it is straightforward to check that
these conservation laws lead to the  following $H^1(\R)$ a priori
estimate for regular solutions of the $(BO)$ equation,
\begin{equation}
  \forall t\geq 0 \; , \; \| u(t)\|_{H^1}  \leq C \; \| u_0\|_{H^1}  \;\; .
  \label{H1conservation}
\end{equation}
Thus, any local well possedness result in $H^1(\R)$ for the $(BO)$
equation can be extend to a global one.

In our knowledge the first results concerning the well possedness of $(BO)$ in the Sobolev spaces $H^s(\R)$ have been obtained in \cite{Saut}
where the global well posedness is proved in $H^3(\R)$. This was improved later to a global well posedness result in $H^s(\R)$, $s>3/2$ in \cite{ABFS},
 \cite{Io} and next  in $H^{3/2}(\R)$, see \cite{Po}. Then, by means of some dispersive estimates for
 the non homogeneous linear Benjamin-Ono equation, $(BO)$ has been proved
 to be locally well posed in $H^s(\R)$, $s>5/4$ in \cite{KT1} and next in $H^s(\R)$, $s>9/8$, \cite{KK}.
 Recently T. Tao  \cite{Tao} get the global
 well posedness in $H^1(\R)$ by using a gauge transformation together with the
  well known Strichartz estimate
\begin{equation}
  \|   V(t)\varphi\|_{L^4_{t,x}} \leq C\, \| \varphi \|_{L^2} \; ,
  \label{Strich}
\end{equation}
where $V(\cdot)$ denotes the free linear group associated to the
linear Benjamin-Ono equation. Up to now the best result concerning
this problem is due to A.D. Ionescu and C.E. Kenig who obtained very
recently the global well posedness of $(BO)$ in $L^2(\R)$,
\cite{IK}.

It is worth noticing that all these recent results have been obtained by coupling compactness methods together with "smoothing" estimates for $V(\cdot)$.
 It is also important to notice that, for all $s \in \R$,  the flow map is not of class $C^2$ from $ H^s(\R)$ to $ C([0,T], H^s(\R))$, \cite{MST}.
 Actually  it has recently been proved in \cite{KT2} that, for all $s>0$,  the flow map is not even uniformly
 continuous on bounded sets of $H^s(\R)$. This is mainly due to some bad interactions between low and high frequencies in the nonlinear
 term $uu_x$ as pointed out in \cite{MST} and \cite{KT2} and this explains why contraction methods can not be used to solve $(BO)$ in $H^s(\R)$.

Concerning the periodic case, the global well posedness in
$H^s(\T)$, $s>3/2$ is derived in \cite{ABFS}. It is worth noticing
that the proof did  not use the smoothing properties of $V(\cdot)$.
Recall that  in the periodic case both the dispersive estimates
\begin{equation}
  \| V(t)\varphi \|_{L^\infty_x } \lesssim t^{-1/2}  \, \| \varphi \|_{L^1} \; ,
  \label{dispers}
\end{equation}
and the sharp Kato smoothing effect
\begin{equation}
 \| D_x^{1/2} V(t) \varphi \|_{L^{\infty}_xL^2_t}\lesssim \| \varphi \|_{L^2}
  \label{Katosmooth}
\end{equation}
fail. This probabely explains why there is no result concerning the
periodic Cauchy problem in $H^s(\T)$, $s\leq 3/2$, for the $(BO)$
equation. In this paper, following the work of T. Tao \cite{Tao}, we
use a gauge transformation together with the periodic estimate
(\ref{Strich}) proved in \cite{Bo1} to obtain the following result
(see subsection 1.2 below for the definition of the space $X^1_T$) :
\begin{theorem}
For all $u_0\in H^1(\T)$ and all $ T> 0 $, there exists a unique
global solution $u$ of the Benjamin-Ono equation in
$$
X^1_T  \, \cap \, C_b(\R , H^1(\T)) \;.
$$
Moreover, the flow-map is continuous  from $H^1(\T)$ to $C([0,T] ,
H^1(\T))$ and, for all $\gamma \in \R$, is Lipschitz on every bounded subset of 
$\Gamma_\gamma$ where
$$
  \Gamma_\gamma =\{ f\in H^1(\T) \, , \, \into f(x) \, dx=\gamma
  \}\, .
$$
\end{theorem}

\begin{rem}
To prove the above Lipschitz property and the uniqueness part of
Theorem 1.1, we will use in a crucial way that, in sharp contrast
with the non periodic case, the gauge transformation is in fact a
Lipschitz map from the set of  $L^2$  functions with zero mean value
on the circle into $L^\infty$. This also avoid to consider some
frequency enveloppe considerations as done in \cite{Tao}.
\end{rem}

When $k \geq 2$, $(GBO)$ is no more a totally integrable system and there is no conservation law at the level $H^1$.
Nonetheless we still have the three following quantities conserved by the flow,
$$
 I(u)=\int \, u(t,x) \, dx \,,\quad M(u)=\int \, u^2(t,x) \, dx \,,
$$
and
$$
  \tilde{E}(u)=\int \, \Bigl( \frac{1}{2} \, |D_x^{1/2}u(t,x)|^2\mp \frac{1}{(k+1)(k+2)}\, u(t,x)^{k+2}\Bigr) \, dx \;\; (energy) \;\;.
$$

In \cite{MR2}, in the case of the whole real line, by means of a
gauge transformation together with some linear estimates, we proved
the local well posedness of  $(GBO)$ on the line in $H^s(\R)$,
$s\geq 1/2$ for $k\geq 5$,  in $H^s(\R)$, $s> 1/2$ for $k=2,4$ and
in $H^s(\R)$, $s\geq 3/4$ for $k=3$ \footnote{See also \cite{MR1}
where optimal results are obtained for $(GBO)$ by contraction
methods in the particular context of small initial data.}. In all
those cases we also obtained that, in a sharp contrast with the
$(BO)$ equation, the flow map is  lipschitz on bounded set of
$H^s(\R)$ which has to be viewed as a stability result for the
$(GBO)$ equation when $k\neq  1$.

In the periodic case, using again the proofs given in \cite{Io}, it
is straightforward to derive the local well posedness of $(GBO)$ in
$H^s(\T)$ when $s>3/2$. On the contrary, up to our knowledge, there
is no available result on this problem when $s \leq 3/2$. As for the
$(BO)$ equation, this is probably due to the failure of the
dispersive estimate (\ref{dispers}) and the sharp Kato smoothing
effect (\ref{Katosmooth}). Again, using a gauge transformation and
the periodic estimate (\ref{Strich}) we prove the following result,
\begin{theorem} Let $k\geq 2 $ be an integer.
For all $u_0\in H^1(\T)$
there exists $ T = T(\|u_0\|_{H^1})>0 $ and a unique solution $u$ of (GBO) in
$$
X^1_T  \, \cap \, C([0,T] ,  H^1(\T)) \;.
$$
Moreover, the flow-map is continuous from $H^1(\T)$ to
$C([0,T] , H^1(\T))$.
\end{theorem}

%%%%%%%%%%%%%%%%%%%%%%%%%%%%%%%%%%%%%%%%%%%%%%%%%%%%%%%%%%%%%%%%%%%%%%%%%%%%%%%%%%%%%%%%%%%%%%%%%%%%%%%
\subsection{Notations}
In the sequel $C$ denotes a positive constant which may differ at each appearance.
 When writing $x\lesssim y$ (for $x$ and $y$ two nonnegative real numbers), we mean that there exists $C_1$ a positive constant (which does not depend of $x$ and $y$) such
that $ x\leq C_1 y$.

We will use the space-time Lebesgues spaces  $L^p_TL^r_\lambda$  of
the $2\pi \lambda$-periodic function (in space) endowed with the
norms

$$
    \|f(t,x)\|_{L^p_TL^r_\lambda}=\| \, \| f(t,\cdot)\|_{L^r([0, 2\pi \lambda])}\|_{L^p([-T,+T])} \;.
$$
When $p=r$, we rather use the notation $L^r_{T,\lambda}=L^r_TL^r_\lambda$.\\
We will also need the functional spaces $X^0_T$, $X^1_T$ and $X^2_T$
respectively defined trough the norms
$$
\|u\|_{X^0_\lambda}=\| u \|_{L^\infty_T L^2_\lambda} +\| u \|_{L^4_{T,\lambda}}  \, ,
$$
$$
\|u\|_{X^1_\lambda}=\|u\|_{X^0_\lambda} + \| u_x \|_{L^\infty_T L^2_\lambda} +\| u_x \|_{L^4_{T,\lambda}}  \, ,
$$
and
$$
\|u\|_{X^2_\lambda}=\|u\|_{X^1_\lambda} + \| u_{xx} \|_{L^\infty_T L^2_\lambda} +\| u_{xx} \|_{L^4_{T,\lambda}}  \, .
$$
In a standard way $H^s_\lambda$ denotes the space of $2\pi
\lambda$-periodic functions such that
$$
  \|f\|_{H^s_\lambda}=  \Bigl(  \sum_{q \in {\Z /\lambda }} \, (1+q^2)^{s}\, | C_q(f)|^2 \Bigr)^{1/2} < +\infty
  \;,
$$
where, for $q\in {\lambda^{-1}\Z  }$,
$$
C_q(f)=\frac{1}{2\pi \lambda } \int_{0}^{2\pi \lambda}
f(t)e^{-iqt}\, dt \;.
$$
 Also, for a $2\pi \lambda$-periodic
function $f$ we respectively defined the projection operators $P_+$,
$P_-$, $P_k$ and $P_{>k}$ by
$$
P_+(f)=\sum_{q \in {\Z_{+}^{*} / \lambda}} \, C_q(f) \, e^{iqx} \; ,
\, P_-(f)=\sum_{q \in {\Z_{-}^{*} / \lambda}} \, C_q(f) \, e^{iqx}
\; ,
$$
$$
P_k(f)=\sum_{q \in {\Z / \lambda}, \, |q|\leq k} \, C_q(f) \,
e^{iqx} \quad \mbox{and} \quad P_{>k}(f)=\sum_{q \in {\Z / \lambda},
\, q > k} \, C_q(f) \, e^{iqx} \; .
$$

\section{Linear estimates}
Let us first recall the following  estimate established by Bourgain
\cite{Bo1} (it was proven for the  Schr\"odinger group but the
adaptation for the Benjamin-Ono group is straightforward).
\begin{equation}
\|V(t) \varphi \|_{L^4_{1} L^4_1} \lesssim
\|\varphi\|_{L^2_1} \quad . \label{Strich1}
\end{equation}
Of course \re{Strich1} still holds for any period $ \lambda \sim 1 $. On the other hand,  by a scaling argument, it is clear that  pushing
$ \lambda $ to $ +\infty $, a factor
 $ \lambda^{1/4} $ will appear at the right-hand side of \re{Strich1}. Since  to solve the Cauchy problem
  for large initial data, we  will use a rescaling argument that will make us work with a large period, the estimate
   \re{Strich1} will not be sufficient.  We will rely instead on the following improved Zygmund estimate also shown in \cite{Bo1} :
 \begin{equation}
\|u\|_{L^4_1 L^4_1} \lesssim \|u\|_{X^{3/8,0}_1} \label{estBourgain} \quad .
\end{equation}
Here the $ X^{s,b}_1 $ are the function spaces introduced by Bourgain
to solve the Cauchy problem for the periodic Schr\"odinger
equation. Recall that
\begin{equation}
\|u\|_{X^{b,s}_1}=\|U(-t) u(t) \|_{H^{b,s}_1} \label{Xsb} \quad ,
\end{equation}
where $ U(t) $ is the Schr\"odinger one parameter linear group.
Again, by separating the positive and the negative frequencies of $u
$ and using \re{estBourgain}, it is obvious to see  that
\re{estBourgain} still holds when replacing $ U(\cdot) $ by $
V(\cdot) $ in the definition of $ X^{s,b} $(see  \re{Xsb}). With
\re{estBourgain} in hand, we establish now the following lemma that
gives a periodic estimate which turns out to be uniform with respect
to large periods $\lambda$.
\begin{lem}
There exists a constant $ C> 0 $ such that $ \forall \lambda\ge 1
 $, $ \forall \varphi\in L^2_\lambda $,
\begin{equation}
\|V(t) \varphi \|_{L^4_{[0,1]} L^4_\lambda} \le C
\|\varphi\|_{L^2_\lambda} \quad . \label{Strich}
\end{equation}
\end{lem}
{Proof .} We first take $ \lambda=1 $. Let $ \psi $ be a
$C^\infty $-function such that $ 0\le \psi\le 1 $, $ \psi \equiv 1
$ on $ [-1/4,1/4] $ and $ \supp \psi \subset [-1/2,1/2] $. We set
$ \psi_T(\cdot)=\psi(\cdot/T) $.
Of course, $ \psi_T \equiv 1 $ on $[-\frac{T}{4},\frac{T}{4}] $ and $ \supp \psi_T \subset [-\frac{T}{2},\frac{T}{2}] $.
 Moreover, for $0<T\le 1 $, $ \psi_T $ can be extended outside $ [-1/2,1/2] $ to a
1-periodic function. Therefore, by \re{estBourgain},
\begin{eqnarray}
\|V(t) \varphi  \|_{L^4_{[-\frac{T}{4},\frac{T}{4}]} L^4_1}
 & \lesssim &\|\psi_T \, V(t) \varphi  \|_{L^4_{1} L^4_1}=
\| V(t) (\psi_T\varphi)  \|_{L^4_{1} L^4_1}\nonumber \\
& \lesssim & \|\psi_T\varphi\|_{X^{3/8,0}_1} \nonumber \\
& \lesssim & \|\psi_T \|_{H^{3/8}_1} \|\varphi\|_{L^2_1}
 \lesssim  T^{1/8} \|\varphi\|_{L^2_1} \quad . \label{k1}
\end{eqnarray}
Where in the last step we use that
$$
\|\psi_T \|_{H^{3/8}_1}\lesssim  \|\psi_T \|_{L^2_1}+\|\psi_T
\|_{\dot{H}^{3/8}_1} \lesssim T^{1/8} \quad .
$$
Now for $ \varphi_\lambda \in L^2_\lambda $, we define
 $ \varphi(\cdot)=\lambda \varphi_{\lambda}(\lambda\cdot)\in L^2_1 $.
 Setting $ u_\lambda(t,x)=(V(t)\varphi_\lambda)(x) $ and
 $ u(t,x)=(V(t) \varphi)(x) $, we easily check that
 $ u(t,x)=\lambda u_\lambda(\lambda^2 t, \lambda x) $. Hence,
 taking $ T= 4\lambda^{-2} $ in  \re{k1}, we get
 \begin{eqnarray}
 \|V(t) \varphi_\lambda \|_{L^4_{[0,1]} L^4_\lambda} & = &
  \lambda^{-1/4}  \|V(t) \varphi \|_{L^4_{[0,1/\lambda^2]}
  L^4_1} \nonumber\\
  & \lesssim & \lambda^{-1/2} \|\varphi\|_{L^2_1} \nonumber \\
  & \lesssim & \|\varphi_\lambda\|_{L^2_\lambda} \quad .
 \end{eqnarray}
This completes the proof of the lemma.

\section{Proof of Theorem 1}
We start by the proof of the local well-posedness result for the
Benjamin-Ono equation which is  much simpler than the one for
$(GBO)$. This is mainly
due to the three following facts :\\
 $\bullet$ The equation satisfied by the gauge transform is simpler. \\
 $\bullet$ There exists a global existence result for smooth initial data.\\
 $\bullet$  The  $ L^2 $-norm is surcritical for $(BO$).\\
 \subsection{Gauge transformation and nonlinear estimates}
 Let $ \lambda \geq 1 $ and
 $ u $ be a global $ H^\infty_\lambda $-solution of $(BO)$ with initial data $ u_0 $.
 In the sequel,  we  assume that $u_0 $ has zero mean value.
  Otherwise we do  the change of unknown
 \begin{equation}
 v(t,x)=u(t,x-t \into u) -\into u \label{chgtvar}\quad  .
 \end{equation}
 Since $ {\displaystyle \into u }$ is preserved by the flow, it is straightforward to see that $ v $ satisfies $(BO)$ with $ {\displaystyle v_0=u_0-\into u_0}$ as
 initial data and that ${\displaystyle \into v=0}$.  Hence we are reduced to the case of a zero mean-value initial data and thus to the case of zero mean-value solutions. Also, changing $u$ to $u/2$, we can always assume that $u$
satisfies the equation,
$$
u_t+{\cal{H}}u_{xx}=2u u_x\; .
$$
 We  define  $ F=\partial_x^{-1} u $ which is the periodic, zero mean value, primitive of $ u $ and
following T. Tao \cite{Tao}, we introduce the gauge transform
$$
W=P_+(e^{-iF}) \quad ,
$$
and we consider
$$
  w=W_x=-iP_+(e^{-iF} F_x) =-iP_+(e^{-iF}
u) \;.
$$
Following the calculations performed in subsection 4.1 (for $k=1$) and noticing that
$$
P_+(e^{-iF}P_-(u_{xx}))-i P_+(u e^{-iF}P_-(u_x))=\partial_xP_+\Bigl(
e^{-iF}P_-(F_{xx}) \Bigr) \, ,
$$
and that
$$
\partial_x^{-1}\, \partial_x(u^2)=u^2-P_0(u^2) \, ,
$$
we then obtain that $w$ solves the following Schr\"odinger equation
:
\begin{eqnarray}
  w_t-iw_{xx} &=& -2 \, \partial_x P_+\Bigl( P_-(F_{xx}) e^{-iF}\Bigr) + P_0(u^2) P_+(u e^{-iF}) \nonumber \\
  &=& -2 \partial_x P_+ \Bigl( P_-(u_{x}) e^{-iF} \Bigr)  + P_0(u^2) P_+(u e^{-iF}) \label{A1} \quad .
\end{eqnarray}
Using the Duhamel formulation of this Schr\"{o}dinger equation on
$w$ it follows from \re{Strich} and standard $TT^{*}$ arguments
that
$$
\| w \|_{X^1_{T,\lambda}}\leq \|w(0)\|_{H^1_\lambda }+ \| -2
\partial_x P_+ \Bigl( P_-(u_{x}) e^{-iF} \Bigr)  + P_0(u^2) P_+(u
e^{-iF})\|_{L^1_TH^1_\lambda}
$$
Now, from Lemma 4.1 and the estimate for $|P_0(u^2)|_{L^{\infty}_T}$
derived in subsection 4.2, we obtain that, for $0<T \le 1 $,
\begin{eqnarray}
   \|w\|_{X^1_{T,\lambda}} &\lesssim & \|w(0)\|_{H^1_\lambda} + \|\partial_x P_+( P_-(u_{x}) W )\|_{L^1_T
   H^1_\lambda}
   \nonumber \\
   &+&|P_0(u^2)|_{L^{\infty}_T} \|P_+(u \, e^{-iF})\|_{L^1_TH^1_{\lambda}}
   \nonumber \\
   & \lesssim & \|w(0)\|_{H^1_\lambda}+ T^{1/2}\, \|u_x\|_{L^4_{T,\lambda}} \|J^1_x w\|_{L^4_{T,\lambda}}
   \nonumber \\
   &+&T \|u\|_{X^1_{T,\lambda}}( \|u\|_{X^1_{T,\lambda}}+\|u\|_{X^1_{T,\lambda}}^2) \nonumber \\
\end{eqnarray}
which leads to
\begin{eqnarray}
   \|w\|_{X^1_{T,\lambda}}& \lesssim &  \|u_0\|_{H^1_\lambda}+\|u_0\|_{H^1_\lambda}^2
   \nonumber \\
   &+&   T^{1/2} \|u\|_{X^1_{T,\lambda}} (\|u\|_{X^1_{T,\lambda}}
   + \|u\|_{X^1_{T,\lambda}}^2+\|w\|_{X^1_{T,\lambda}}) \quad . \label{A2}
\end{eqnarray}
On the other hand, we can rewrite $ u $ as
\begin{equation}
 u =  e^{iF} e^{-iF} u  =
e^{iF} P_+(e^{-iF} u) +e^{iF} P_- (e^{-iF} u) \quad , \label{A3}
\end{equation}
and so,
\begin{eqnarray*}
P_{>1} u & = & i P_{>1} \Bigl( e^{iF}w\Bigr) +P_{>1} \Bigl( e^{iF} P_-(e^{-iF} u ) \Bigr) \nonumber \\
& =  &  i P_{>1} \Bigl( e^{iF}w\Bigr) +P_{>1} \Bigl(
P_{> 1}(e^{iF}) {P_-}(e^{-iF} u ) \Bigr) \quad .
\end{eqnarray*}
Hence  from Lemma 4.1 (and since  $ u $ is real-valued),  we infer
that
\begin{eqnarray}
     \|u\|_{X^1_{T,\lambda}} & =  & \|P_1 u\|_{X^1_{T,\lambda}}+ 2\|P_{>1}\,  u \|_{X^1_{T,\lambda}} \nonumber \\
     & \lesssim &  \|P_1 u\|_{X^1_{T,\lambda}}+ \|w\|_{X^1_{T,\lambda}}+ \|u\|_{X^1_{T,\lambda}}\|w\|_{X^0_{T,\lambda}} \nonumber \\
     &+&\|P_{>1} (e^{iF})\|_{L^\infty_{T,\lambda}} \, \|u\|_{X^1_{T,\lambda }}+\|u\|_{X^1_{T,\lambda}}^2 \nonumber \\
    &\lesssim& \|P_1 u\|_{X^1_{T,\lambda}}+ \|w\|_{X^1_{T,\lambda}}+\|u \|_{L^\infty_{T,\lambda}} \, \|u\|_{X^1_{T,\lambda}}+\|u\|_{X^1_{T,\lambda}}^2 \label{A4}
\end{eqnarray}
where in the last step we use that, by  Bernstein's inequality for $
2\pi \lambda $-periodic functions,
$$
\| P_{>1} (e^{iF})\|_{L^\infty_{T,\lambda }} \lesssim \|\partial_x (e^{iF})\|_{L^\infty_{T,\lambda }}
=\|e^{iF} F_x\|_{L^\infty_{T,\lambda }}=\|F_x\|_{L^\infty_{T,\lambda }} =\| u \|_{L^\infty_{T,\lambda }} \, ,
$$
and that, $ \|w\|_{X^0_{T,\lambda}} \le
\|u\|_{X^0_{T,\lambda}} $.
\subsection{Local well-posedness for small data}
We will now prove the local well-posedness result for small initial data. The result for arbitrary large data  will
  follow from scaling arguments.

\subsubsection{Existence}
Let $ u_0\in H^\infty_\lambda $ be a $ 2\pi \lambda$-periodic zero
mean-value function and let us assume that $ \|u_0\|_{H^1_\lambda}
\lesssim \varepsilon^2 $ for some small $0<\varepsilon<\! \! < 1$
depending only on the implicit constant contained in the above
estimates. At this stage, it is worth recalling that these  implicit
constants do not depend on the period $ \lambda$.

Our aim is to show that the emanating solution
$ u\in C(\R;H^\infty_\lambda) $ satisfies  $ \|u\|_{X^1_{1,\lambda }} \lesssim \varepsilon^2 $.\\
First, since the $ L^2 $-norm of $ u $ is  constant along the trajectory, we obtain from Bernstein inequalities that
\begin{equation}
 \|P_1 u \|_{X^1_{1,\lambda}} \lesssim \|u\|_{L^\infty_1 L^2_x} = \|u_0 \|_{ L^2_x} \lesssim \varepsilon^2 \label{A5} \quad .
 \end{equation}
 On the other hand, since $\| w(0) \|_{H^1_\lambda}\lesssim \|u_0\|_{H^1_\lambda}+\|u_0\|_{H^1_\lambda}^2 \lesssim \varepsilon ^2$, by continuity we can  assume that
 $$
\|u\|_{X^1_{T,\lambda }}\lesssim \varepsilon \hspace*{5mm} \mbox{
and } \hspace*{5mm} \|w\|_{X^1_{T,\lambda }}\lesssim \varepsilon \;
,
$$
on some small enough interval $ [0,T] \subset [0,1] $.
 But \re{A2} then gives $ \|w\|_{X^1_{T,\lambda}} \lesssim \varepsilon^2 $ and this last inequality together with \re{A4}-\re{A5} imply now that
  $ \|u\|_{X^1_{T,\lambda}} \lesssim \varepsilon^2 $. In a standard way this proves that $ T $ can be taken equal to $ 1 $ and thus we have
  $$
\|u\|_{X^1_{1,\lambda}} \lesssim \varepsilon^2 \quad .
$$

Consider now $ u_0\in H^1_\lambda $ such that $
\|u_0\|_{H^1_\lambda} \lesssim \varepsilon^2 $. Approximating $ u_0
$ in $ H^1_\lambda $ by a sequence $ \{u_{0,n}\}
 \subset H^\infty_\lambda $, it follows that the sequence of the
  emanating solutions $\{ u_n\} \subset C(\R;H^\infty_\lambda) $ is
  bounded in $ X^1_{1,\lambda} $. We can thus pass to the limit up to a subsequence
   and obtain the existence of a solution $ u\in X^1_{1,\lambda} $ of $(BO)$.

\subsubsection{Continuity, uniqueness and  regularity of the flow map }
As we notice already in the introduction, one of the main differences with
the problem on the real axis is that the gauge transformation is
Lipschitz from the space of  $ L^2 $ functions with zero mean
value on the circle into
 $ L^\infty $. This property turns out to be crucial to get the uniqueness and the continuity of the flow.
We   first prove  that the flow-map is Lipschitz on a small ball of
$ H^1_\lambda $. The continuity of $ t\mapsto u (t) $ in $
H^1_\lambda $ will follow directly.

Let  $ u_1 $ and $ u_2 $ be two solutions of $(BO)$ in $
X^1_{T,\lambda} $ associated with the initial data $ \varphi_1 $ and
$ \varphi_2 $
 in $ H^1_\lambda $.
 We assume that they satisfy
 \begin{equation}
 \|u_i\|_{X^1_{T,\lambda}}\lesssim \varepsilon^2 \; , \; \quad i=1,2, \label{A6}
 \end{equation}
for some $ 0<T\le 1 $ and  where $ \varepsilon $ is taken as
above. We set
 $$
 z=w_1-w_2=-i P_+(e^{-iF_1} u_1)+iP_+ (e^{-iF_2}
  u_2)
 $$
 with $ F_i $ is defined as $ F $ and where $ u $ is replaced by $ u_i $. Obviously, $ z $ satisfies
 \begin{eqnarray}
 z_t-i z_{xx} &=&\partial_x P_+ \Bigl[ P_-(\partial_x u_1 -\partial_x u_2) W_1\Bigr]
                +\partial_x P_+ \Bigl[ P_-(\partial_x u_2)( W_1-W_2)\Bigr] \nonumber \\
              &+& P_0(u_1^2) \, P_+\Bigl((u_1-u_2)e^{-iF_1}\Bigr) + P_0(u_1^2) \, P_+\Bigl(u_2(W_1-W_2)\Bigr) \nonumber \\
              &+& P_0\Bigl(z(u_1+u_2)\Bigr) \, P_+(u_2) \; . \label{A7}
 \end{eqnarray}
 Note that \re{A6} clearly ensures that for $ i=1,2 $,
 \begin{equation}
 \|w_i\|_{X^1_{T,\lambda}} \lesssim
 \|u_i\|_{X^1_{T,\lambda}} (1+ \|u_i\|_{X^1_{T,\lambda}} )\lesssim
 \varepsilon^2 \, . \label{A6b}
 \end{equation}
 As previously, it follows that
 \begin{equation}
 \|z\|_{X^1_{T,\lambda}} \lesssim \|z(0)\|_{H^1_\lambda} + \| A \|_{L^1_TH^1_{\lambda}} \label{A8}
 \end{equation}
 where $A$ denotes the right hand side of \re{A7}.
Note first that,
\begin{eqnarray}
\|z(0)\|_{H^1_\lambda} & \lesssim &
\|\varphi_1-\varphi_2\|_{H^1_\lambda}\Bigl(1+\|\varphi_1\|_{H^1_\lambda}
 +\|\varphi_2\|_{H^1_\lambda}\Bigr)\nonumber\\
 &+& \|e^{-iF_1(0)}-e^{-iF_2(0)}\|_{L^\infty} \|\varphi_2\|_{H^1_\lambda}
 (1+\|\varphi_2\|_{H^1_\lambda})
\label{A9}
\end{eqnarray}
with
\begin{eqnarray}
\|e^{-iF_1(0)}-e^{-iF_2(0)}\|_{L^\infty_{\lambda}} \lesssim   \|\partial_x^{-1}(\varphi_1-\varphi_2 )\|_{L^\infty_\lambda}
 \lesssim    \lambda^{1/2} \, \|\varphi_1-\varphi_2 \|_{L^2_\lambda} \label{A1}
\end{eqnarray}
and so,
\begin{equation}
  \|z(0)\|_{H^1_\lambda}  \lesssim (1+ \lambda ^{1/2} \varepsilon
  ^2) \|\varphi_1-\varphi_2\|_{H^1_\lambda}\, . \label{A1bis}
\end{equation}
We give now an estimate for $\|A \|_{L^1_TH^1_\lambda}$. From Lemma
4.1 we easily obtain that
 \begin{eqnarray}
   \| \partial_x P_+ \Bigl[ P_-(\partial_x u_1 -\partial_x u_2) W_1\Bigr] \|_{L^1_T H^1_\lambda }
   &\lesssim& T^{1/2} \, \|u_1-u_2\|_{X^1_{T, \lambda }} \, \| w_1\|_{X^1_{T,\lambda }}  \, ,
   \nonumber \\
   \| \partial_x P_+ \Bigl[ P_-(\partial_x u_2)( W_1-W_2)\Bigr] \|_{L^1_T H^1_\lambda }
   &\lesssim& T^{1/2} \, \|u_2\|_{X^1_{T, \lambda }} \, \| z \|_{X^1_{T,\lambda }}   \, ,
   \nonumber  \\
   \| P_0(u_1^2) \, P_+\Bigl[(u_1-u_2)e^{-iF_1}\Bigr] \|_{L^1_T H^1_\lambda }
   &\lesssim &T^{1/2} \, \|u_1-u_2\|_{X^1_{T, \lambda}}(1+\|u_1\|_{X^1_{T, \lambda }})\|u_1\|_{X^1_{T,\lambda}}^2 \, ,
   \nonumber \\
   \| P_0\Bigl(z(u_1+u_2)\Bigr) \, P_+(u_2) \|_{L^1_T H^1_\lambda }
   &\lesssim& T^{1/2} \, \|z\|_{X^1_{T,\lambda}} \|u_2\|_{X^1_{T, \lambda}}
   (\|u_1\|_{X^1_{T, \lambda }}+\|u_2\|_{X^1_{T, \lambda }})\, , \nonumber
   \label{A10bis}
 \end{eqnarray}
 and proceeding as for $\|z(0)\|_{H^1_\lambda}$,
 \begin{eqnarray}
   \| P_0(u_1^2) \, P_+(u_2(W_1-W_2)) \|_{L^1_T H^1_\lambda }
   &\lesssim& T^{1/2} \, \|u_1\|_{X^1_{T,\lambda}}^2
   \|u_2\|_{X^1_{T,\lambda}}\|W_1-W_2\|_{L^{\infty}_{T,x}} \nonumber \\
   &+& T^{1/2} \, \|u_1\|_{X^1_{T,\lambda}}^2 \|u_2\|_{X^1_{T,\lambda}} \|z\|_{X^1_{T,\lambda}}\nonumber \\
   &\lesssim& T^{1/2} \lambda^{1/2} \, \|u_1\|_{X^1_{T,\lambda}}^2
   \|u_2\|_{X^1_{T,\lambda}}\|u_1-u_2\|_{L^{\infty}_TL^2_x} \nonumber \\
   &+& T^{1/2} \, \|u_1\|_{X^1_{T,\lambda}}^2 \|u_2\|_{X^1_{T,\lambda}} \|z\|_{X^1_{T,\lambda}}\, .\label{A11bis}
 \end{eqnarray}
 Hence gathering \re{A6b}, \re{A1bis} and the previous estimates we infer that,
\begin{eqnarray}
  \|z\|_{X^1_{T,\lambda}} &\lesssim&
  (1+\varepsilon ^2 \lambda^{1/2})\|\varphi_1-\varphi_2\|_{H^1_\lambda }\nonumber \\
   &+& \varepsilon ^2 T^{1/2}\Bigl( \|z\|_{X^1_{T,\lambda}}+(1+\lambda
  ^2)\|u_1-u_2\|_{X^1_{T,\lambda }}\Bigr) \label{A12bis} \, .
\end{eqnarray}
On the other hand, we have
\begin{eqnarray*}
u_1-u_2 & =& \partial_x F_1 -\partial_x F_2 \nonumber \\
 & = & i e^{iF_1} \Bigl[ z +\partial_x P_- \Bigl( e^{-iF_1}  -e^{-iF_2}
 \Bigr)\Bigr] +i ( e^{iF_1} -e^{iF_2}) \Bigl(w_2
  +  \partial_x P_- (e^{-iF_2} )\Bigr) \label{A10}
\end{eqnarray*}
and thus
 \begin{eqnarray}
P_{+} (u_1-u_2)
 & = & i P_{+}(e^{iF_1}z)
 +i P_{+}\Bigl[  e^{iF_1}\partial_x P_-\Bigl( e^{-iF_1} -e^{-iF_2} \Bigr)\Bigr] \nonumber \\
& &   +i  P_{+}\Bigl[( e^{iF_1} -e^{iF_2})   w_2\Bigr]
  + i P_{+}\Bigl[  (e^{iF_1} -e^{iF_2}) \partial_x P_- (e^{-iF_2} ) \Bigr] \label{A11}  \quad . \nonumber
\end{eqnarray}
 Therefore, as in \re{A4},
\begin{eqnarray}
\| u_1-u_2\|_{X^1_{T,\lambda}}  & \lesssim
&\|z\|_{X^1_{T,\lambda}}(1+\|u_1\|_{X^1_{T,\lambda}} )
+\|u_1-u_2\|_{X^1_{T,\lambda}} \|u_2\|_{X^1_{T,\lambda}}\nonumber \\
 & & +\|e^{iF_1}-e^{iF_2}\|_{L^\infty_{T,\lambda}}
 \Bigl( \|u_1\|_{X^1_{T,\lambda}}+\|u_2\|_{X^1_{T,\lambda}}+
 \|u_1\|_{X^1_{T,\lambda}}\|u_2\|_{X^1_{T,\lambda}}\Bigr)
 \nonumber \\
 & & + \|w_2\|_{X^1_{T,\lambda}}\Bigl( \|e^{iF_1}-e^{iF_2}\|_{L^\infty_{T,\lambda}}(1+\|u_1\|_{X^1_{T,\lambda}})
 +  \|u_1-u_2\|_{X^1_{T,\lambda}}\Bigr)  \quad . \nonumber \label{A12}
\end{eqnarray}
But, proceeding as in \re{A1}, we see that
\begin{equation}
\|e^{-iF_1}-e^{-iF_2}\|_{L^\infty_{T,\lambda}}
 \lesssim \lambda^{1/2} \|u_1-u_2\|_{L^\infty_T L^2_{\lambda}} \quad . \label{A12b}
\end{equation}
Writing now the equation satisfied by  $ u_1-u_2 $ and using it's
Duhamel formulation together with $TT^*$ arguments and Sobolev
inequalities we easily obtain that
  for $ 0<T\le 1 $,
\begin{eqnarray}
\|u_1-u_2\|_{L^\infty_T L^2_\lambda}
  \lesssim  \|\varphi_1 -\varphi_2\|_{L^2_\lambda} +T^{1/2}
\|u_1-u_2\|_{X^1_{T,\lambda}}
 (\|u_1\|_{X^1_{T,\lambda}}+\|u_2\|_{X^1_{T,\lambda}}) \quad . \label{A12c}
\end{eqnarray}
 Gathering these estimates and recalling
\re{A6}, \re{A6b} we finally obtain
\begin{eqnarray}
\|u_1-u_2\|_{X^1_{T,\lambda}} & \lesssim & (1+
\lambda^{1/2}\varepsilon^2) \|\varphi_1-\varphi_2\|_{H^1_\lambda}
 \nonumber \\
 & &+ \varepsilon^2 \|z\|_{X^1_{T,\lambda}}+\varepsilon^2 \|u_1-u_2\|_{X^1_{T,\lambda}}
 + \lambda^{1/2} \varepsilon^2 \|u_1-u_2\|_{L^\infty_T L^2_\lambda} \nonumber \\
  &\lesssim & (1+
\lambda^{1/2}\varepsilon^2) \|\varphi_1-\varphi_2\|_{H^1_\lambda}+\varepsilon^2 \|u_1-u_2\|_{X^1_{T,\lambda}} \nonumber \\
  & & + T^{1/2} \, \varepsilon^4 \lambda^{1/2}  \|u_1-u_2\|_{X^1_{T,\lambda}} \quad .
\label{A13}
\end{eqnarray}
Hence, for $ 0<T\le  T_\lambda \sim   \lambda^{-1} $, we get
\begin{equation}
\|u_1-u_2\|_{X^1_{T,\lambda}}  \lesssim  (1+ \varepsilon^2
\lambda^{1/2} ) \|\varphi_1-\varphi_2\|_{H^1_\lambda} \label{A14}
\quad .
\end{equation}
With \re{A14} in hand, we observe that the approximative sequence
$ u^n $ constructed above is a Cauchy sequence
 in  $ C([0,T_\lambda]; H^1_\lambda) $ since $  \|u_{n}\|_{X^1_{1,\lambda}} \lesssim \varepsilon^2 $
 and so
 $ u_{0,n} $ converges to $ u_0 $ in $H^1_\lambda$.
  Hence, $ u $ belongs to
  $ C([0,T_\lambda]; H^1_\lambda) $. Repeating this argument we  get that actually $ u\in C([0,1]; H^1_\lambda) $.
Moreover, \re{A14} clearly ensures the uniqueness in the considered
class  and
 that  the  flow-map is
Lipschitz from the ball of $ H^1_\lambda $ with radius $
\varepsilon^2 $ into
  $ C([0,1]; H^1_\lambda) $.
  %%%%%%%%%%%%%%%%%%%%%%%%%%%%%%%%%%%%%%%%%%%%%%%%%%%%%%%%%%%%%%
\subsection{The case of arbitrary large initial data}
Here we used the dilation symmetry of the equation to extend the
result for arbitrary large data. First note that if $ u(t,x) $ is a
$2\pi$-periodic solution of $(BO)$ on $[0,T] $ with initial data $
u_0 $ then $ u_\lambda(t,x)=\lambda^{-1}
u(\lambda^{-2}t,\lambda^{-1} x)$ is a $ 2\pi\lambda $-periodic
solution of $(BO)$ on
 $ [0,\lambda^2 T] $
 emanating from $u_{0,\lambda}=\lambda^{-1}\,  u_0(\lambda^{-1} x) $.

Now, let $ u_0\in H^1 $. If $ \|u_0\|_{H^1}\le \varepsilon^2 $ we
are in the small initial data case. Otherwise, we set
$$
\lambda=\varepsilon^{-4} \|u_0\|_{H^1}^2\ge 1
$$
so that $ u_{0, \lambda} $ satisfies
$$
\|u_{0,\lambda}\|_{H^1_\lambda} \lesssim \varepsilon^2 \quad .
$$
Hence we are reduce to the case of small initial data. Therefore, there exists a unique
 local solution $ u_\lambda\in C([0,1];H^1_\lambda)\cap X^1_{1,\lambda} $ of $(BO)$    emanating from $ u_{0,\lambda} $. This proves the existence
 and uniqueness in $ C([0,T]; H^1)\cap X^1_{T,1} $ of the solution $ u $ emanating from $ u_0 $ with $ T\sim
 \|u_0\|_{H^1}^{-4} $. The fact that the flow-map is Lipschitz on every bounded
 set of $ H^1 $ follows as well.\vspace*{3mm}

 Finally, note that the change of unknown \re{chgtvar} preserves the continuity of the solution and the continuity of the flow-map in $ H^1(\T) $.
 Moreover, the Lipschitz property (on bounded sets) of the flow-map is also preserved on the hyperplans of $ H^1(\T) $ with fixed  mean-value.

\section{Gauge transform for $(GBO)$ and nonlinear estimates}
Let us now begin the proof of Theorem 2. As for the $(BO)$ equation
we have to perform a gauge transformation to obtain suitable
estimates in $X^1_{T,\lambda}$ for regular solutions of $(GBO)$.
\subsection{The gauge transformation}
For $u$ a smooth $2\pi \lambda$-periodic solution of $(GBO)$ we
consider $v$ defined as
\begin{equation}
 v(t,x)=2^{1/k}\, u( t, x+\int_0^t \into u^k) \,,
 \label{chgtvarGBO}
\end{equation}
which satisfies
$$
 \left\{\begin{array}{lll}
\partial_t v +{\cal{H}}\partial^2_x v = 2\,  M(v^k)  \partial_xv \,  \, ,\; (t,x)\in \R\times \T\;,\\
v(0,x)=u_0(x) \;,
\end{array}\right.
$$
where  ${\displaystyle M(g)=g - \into g}$ (such a manipulation is used in \cite{CKSTT} for the generalized Korteweg-de-Vries equations on the torus). In the same spirit as in \cite{HaOz} and \cite{MR2}, define $w$ the gauge transform of $u$ by
\begin{equation}
w=P_+( e^{-iF}v ) \, , \; F(t,x)=\sum_{q \in {\Z^* / \lambda}}
C_q(M(v^k)) \frac{e^{iqx}}{iq}=\partial_x^{-1}(M(v^k)) \;.
\label{defw}
\end{equation}
In the sequel our aim is to  derive a suitable equation satisfied by $w$. Noticing that
$$
\left\{ \begin{array}{lll}
   w_t=P_+[ e^{-iF}(-iF_tv+v_t)] \, ,  \\
   w_{xx}=P_+[e^{-iF}(-2iv_xF_x+v_{xx}-F_x^2v-iF_{xx}v)]  \,,
\end{array}\right.
$$
we obtain that $w$ solves the semilinear Schr\"odinger equation
\begin{equation}
 w_t-iw_{xx}=P_+[e^{-iF}((v_t-iv_{xx}-2F_xv_x)+(-F_{xx}v+iF_x^2v)-iF_tv )]\;.
 \label{eqw}
\end{equation}
We compute now the three terms appearing in the right hand side of
(\ref{eqw}). First we have,
\begin{eqnarray}
  A&=&v_t-iv_{xx}-2F_xv \nonumber \\
  &=&v_t+Hv_{xx}-2iP_-(v_{xx})-2M(v^k)v_x \nonumber \\
  &=&-2iP_-(v_{xx})\,.
\label{defA}
\end{eqnarray}
Next we have,
\begin{eqnarray}
  B&=&-F_{xx}+iF_x^2v \nonumber \\
   &=&-[M(v^k)]_xv+i[M(v^k)]^2v \nonumber \\
   &=&B_1+B_2 \,.
  \label{defB}
\end{eqnarray}
On the other  hand $C=-iv F_t$ with,
\begin{eqnarray}
    F_t  &=&\sum_{q \in {\Z^* / \lambda}}C_q(M(v^k)_t) \frac{e^{iqx}}{iq}\nonumber \\
   &=& k \sum_{q \in {\Z^* / \lambda}}C_q(v^{k-1}v_t) \frac{e^{iqx}}{iq}  \nonumber \\
   &=& - k \sum_{q \in {\Z^* / \lambda}}C_q(v^{k-1}Hv_{xx}) \frac{e^{iqx}}{iq}
   +    2k \sum_{ q \in {\Z^* / \lambda}}C_q(v^{k-1}M(v^k)v_x) \frac{e^{iqx}}{iq} \nonumber \\
   &=& C_1+C_2 \,,
  \label{defC}
\end{eqnarray}
with,
\begin{eqnarray}
  C_1 &=&-k \sum_{q \in {\Z^* / \lambda}}C_q(v^{k-1}Hv_{xx}) \frac{e^{iqx}}{iq}  \nonumber \\
      &=& k \sum_{q \in {\Z^* / \lambda}} [ \sum_{r,q-r \in {\Z^* / \lambda}} (q-r)^2C_r(v^{k-1})C_{q-r}(H v) ]  \frac{e^{iqx}}{iq} \nonumber \\
      &=& -k \sum_{q \in {\Z^* / \lambda}} [ \sum_{r,q-r \in {\Z^* / \lambda}}C_r(v^{k-1})C_{q-r}(Hv_x) ]  e^{iqx}  \nonumber \\
      &+& k(k-1) \sum_{q \in {\Z^* / \lambda}} [ \sum_{r,q-r \in {\Z^* / \lambda}} C_r(v^{k-2}v_x)C_{q-r}(Hv_x) ]  \frac{e^{iqx}}{iq}  \nonumber \\
      &=& -k \sum_{q \in {\Z^* / \lambda}} C_n(v^{k-1}Hv_x)   e^{inx}
      + k(k-1) \sum_{q \in {\Z^* / \lambda}}  C_q(v^{k-2}v_xHv_x)   \frac{e^{iqx}}{iq}  \nonumber \\
      &=&-k M(v^{k-1}Hv_x)+k(k-1) \sum_{ q \in {\Z^* / \lambda}}  C_q( v^{k-2} v_x Hv_x)   \frac{e^{iqx}}{iq} \nonumber \\
       &=&C_{1,1}+C_{1,2} \,.
      \label{defC1}
\end{eqnarray}
and,
\begin{eqnarray}
  C_2 &=&2k \sum_{ q \in {\Z^* / \lambda}}C_q(v^{k-1}M(v^k)v_x) \frac{e^{iqx}}{iq} \nonumber \\
      &=&2  \sum_{ q \in {\Z^* / \lambda}}C_q(M(v^k)M(v^{k})_x) \frac{e^{iqx}}{iq} \nonumber \\
      &=&  \sum_{ q \in {\Z^* / \lambda}}C_q((M(v^k)^2)_x) \frac{e^{iqx}}{iq} \nonumber \\
      &=&M[M(v^k)^2] \, .
     \label{defC2}
\end{eqnarray}
Gathering \re{defA}-\re{defC2}and noticing that
\begin{eqnarray}
   B_1-iv \, C_{1,2}&=& M(v^k)_xv+ikvM(v^{k-1}Hv_x) \nonumber \\
   &=&k v M(v^{k-1}P_+v_x) -k v M(v^{k-1}P_-v_x)-kv^kv_x \nonumber  \\
   &=& k v M(v^{k-1}P_+v_x) -k v M(v^{k-1}P_-v_x)-k v M(v^{k-1} v_x)  \nonumber \\
   &=& -2 k v M(v^{k-1}v_x) \; ,
   \label{defB1}
\end{eqnarray}
and that
\begin{eqnarray}
   B_2-iv \, C_2&=&-iv \, [ M(M(v^k)^2) -M(v^k)^2] \nonumber  \\
   &=& -iv \, P_0(M(v^k)^2) \, ,
  \label{defB2}
\end{eqnarray}
 we infer that,
\begin{eqnarray}
   w_t-iw_{xx} &=& P_+(e^{-iF}v \, P_0(M(v^k)^2)  \nonumber \\
               &-& 2iP_+(e^{_iF}P_-v_{xx}) \nonumber  \\
               &-&2k P_+(e^{-iF} v \, M(v^{k-1}P_-v_x)) \nonumber \\
               &-&k(k-1) P_+\Bigl( e^{-iF} v \sum_{q \in {\Z^* / \lambda}}\frac{C_q(v^{k-2} v_x Hv_x)}{iq} \, e^{iqx} \Bigr) \nonumber \\
               &:=& a+b+c+d \, .
               \label{eqonw}
\end{eqnarray}
which gives us the Schr\"odinger equation satisfies by $w$.

\subsection{Nonlinear estimates on the gauge transform}
In this subsection we derived $X^1_{T,\lambda}$ and
$X^2_{T,\lambda}$-estimates for $w$ the gauge transform of $v$ a
$2\pi \lambda$-periodic solution of $(GBO)$. Our aim is to prove the
following estimate,

\begin{pro}
Let $v$ a $X^1_{T,\lambda}$-solution of $(GBO)$ and let us consider
$w=P_+(e^{-iF}v)$ . Then  for $0\leq T \leq 1$ and for $\lambda >1$,
\begin{equation}
\| w\|_{X^1_{T,\lambda}} \lesssim \| w_0\|_{H^1_{\lambda}} \, + \,
T^{1/4} \, (\| v\|_{X^1_{T,\lambda}}^{k+1}+\|
v\|_{X^1_{T,\lambda}}^{2k+1}+\| v\|_{X^1_{T,\lambda}}^{3k+1} ) \; .
\label{estprop}
\end{equation}
Moreover there exists a polynomial function $Q$ such that,
\begin{equation}
\| w\|_{X^2_{T,\lambda}} \lesssim \| w_0\|_{H^1_{\lambda}} \, + \,
T^{1/4} \, Q(\| v\|_{X^1_{T,\lambda}}) \, \| v\|_{X^2_{T,\lambda}}
\; . \label{estbisprop}
\end{equation}
\end{pro}

\noindent {\it Proof.} We will only prove \re{estprop} since
\re{estbisprop} can be derived in exactly the same way (the linear
dependence in the strong norm $\| u\|_{X^2_{T,\lambda}}$ of  the
right hand side of \re{estbisprop} is standard).  Recall first that
$w$ solves
$$
  w_t-iw_{xx}=f=a+b+c+d \, ,
$$
where $a$, $b$, $c$ and $d$ are defined in \re{eqonw}. Hence from
the periodic estimate
\begin{equation}
\|w \|_{L^{\infty}_tL^2_x}+\|w\|_{L^4_{t,x}}\leq \|f\|_{L^1_tL^2_x}
 \label{Stic1}
\end{equation}
we infer that
\begin{eqnarray}
   \| w\|_{X^1_{T,\lambda}}& \lesssim & \|a\|_{L^1_tL^2_x}+\|\partial_x a\|_{L^1_tL^2_x}+ ... +\|d\|_{L^1_tL^2_x}+\|\partial_x
d\|_{L^1_tL^2_x}\, .
\end{eqnarray}
Now we recall the following periodic version of  of Lemma 3.4 in \cite{MR1} which allows to share derivatives when
estimating terms like $D_x^{\alpha} P_+[fP_-(D^{\beta}_xg)]$.
%%%%%%%%%%%%%%%%%%%%%%%%%%%%       LEMME DERIVEES %%%%%%%%%%%%%%%%%%%%%%%%%%%%%%%%%%%%%%%%%%%%%%
\begin{lem}
Let $\alpha>0$, $\beta\geq 0$, $1<p< +\infty$ and $1\leq q \leq  +\infty$. Then
\begin{equation}
  \|D^{\alpha}_xP_+[fP_-(D^{\beta}_xg]\|_{L^q_qL^p_x}\lesssim \|D^{\gamma_1}_xf\|_{L^{q_1}_tL^{p_1}_x}
  \|D^{\gamma_2}g\|_{L^{q_2}_tL^{p_2}_x}\, ,
\end{equation}
with $1<p_i,q_i<+\infty$, $1/{p_1}+1/{p_2}=1/p$, $1/{q_1}+1/{q_2}=1/q$ and $\gamma_1\geq \alpha$, $\gamma_1 +
\gamma_2=\alpha+\beta$.
\label{lemderiv}
\end{lem}
%%%%%%%%%%%%%%%%%%%%%%%%%%%%%%%%%%%%%%%%%%%%%%%%%%%%%%%%%%%%%%%%%%%%%%%

\noindent {\bf .} Estimates for $a=P_0[M(v^k)^2] \, P_+(e^{-iF}v) \,
$. We have,
\begin{eqnarray}
   \|a\|_{L^1_TL^2_x} & = & |P_0[M(v^k)^2] | \,  \| v \|_{L^1_TL^2_x} \nonumber \\
               & \lesssim & T \, \| M(v^k) \|_{L^{\infty}_TL^2_x}^2 \,  \|v\|_{L^{\infty}_TL^2_x} \nonumber \\
               & \lesssim & T \, \| v\|_{L^{\infty}_TL^{2k}_x}^{2k} \, \|v\|_{L^{\infty}_TL^2_x} \nonumber \\
               & \lesssim & T \, \| v\|_{X^1_{T,\lambda}}^{2k+1} \, .\\
              \nonumber
\end{eqnarray}
Next we have
\begin{eqnarray}
   \|\partial_x a\|_{L^1_tL^2_x} & \lesssim & |P_0[M(v^k)^2] | \,  \| \partial_x P_+(e^{-iF}v )\|_{L^1_tL^2_x} \nonumber \\
               & \lesssim & \| v \|_{X^1_{T,\lambda}}^{2k} \bigl(  \|M(v^k) v\|_{L^1_tL^2_x}+ \|e^{-iF} v_x\|_{L^1_tL^2_x}\bigr) \nonumber \\
               & \lesssim &  \, \| v\|_{X^1_{T,\lambda}}^{2k}
                \bigl( \|M(v^k) \|_{L^\infty _{t,x}} \|v\|_{L^1_tL^2_x}+ \| v_x\|_{L^1_tL^2_x}\bigr) \nonumber \\
               & \lesssim &  T \, \| v\|_{X^1_{T,\lambda}}^{2k}
                \bigl( \| v \|_{X^1_{T,\lambda}}^{k+1} + \|v\|_{X^1_{T,\lambda}} \bigr) \, ,
\end{eqnarray}
where we use that $\|M(v^k) \|_{L^\infty _{t,x}} \lesssim \|M(v^k) \|_{L^\infty _tH^1_x} \lesssim \| v \|_{L^\infty _tH^1_x}^k$.

\noindent {\bf .} Estimates for $b=-2iP_+ [ e^{-iF}P_-(v_{xx}) ]$.
>From Lemma \ref{lemderiv} we have,
\begin{eqnarray}
   \|b\|_{L^1_tL^2_x} & \lesssim & \| \partial_x  e^{-iF} \|_{L^\infty_tL^4_x} \, \| v_x \|_{L^1_tL^4_x} \nonumber \\
               & \lesssim & T^{3/4} \, \| M(v^k) \|_{L^\infty_tL^4_x}  \, \| v_x \|_{L^4_{t,x}} \nonumber \\
               & \lesssim &  T^{3/4} \, \| M(v^k) \|_{L^\infty_tH^1_x} \, \| v_x \|_{L^4_{t,x}} \nonumber \\
               & \lesssim & T^{3/4} \, \| v\|_{X^1_{T,\lambda}}^{k+1}  \, .
\end{eqnarray}
Now, again from Lemma \ref{lemderiv},  we have
\begin{eqnarray}
   \|\partial_x b\|_{L^1_tL^2_x} & \lesssim &   \| \partial_{xx}(e^{-iF})\|_{L^4_{t,x}} \| v_{x}  \|_{L^{4/3}_tL^4_x}   \nonumber \\
                                 & \lesssim & T^{2/3} \, \| v^{k-1} v_{x}  \|_{L^4_{t,x}} \,  \| v_{x}  \|_{L^4_{t,x}}    \nonumber \\
                                 & \lesssim &  T^{2/3} \, \| v^{k-1} \|_{L^\infty_{t,x}}   \,  \|v_{x}  \|_{L^4_{t,x}}^2   \nonumber \\
                                 & \lesssim &  T^{2/3} \, \| v \|_{L^\infty_tH^1_x}^{k-1} \|v_{x}\|_{L^4_{t,x}}     \nonumber \\
                                 & \lesssim &  T^{2/3} \, \| v \|_{X^1_{T,\lambda}}^{k+1}   \, .
\end{eqnarray}

\noindent {\bf .} Estimates for $c=-2kP_+\bigl( e^{-iF}v \,
M[v^{k-1}P_-(v_{x})] \bigr)$.
\begin{eqnarray}
   \|c\|_{L^1_tL^2_x} & \lesssim & \| v \, M[v^{k-1} P_-(v_{x})]\|_{L^1_tL^2_x} \nonumber \\
               & \lesssim & \| v \|_{L^{\infty}_{t,x}} \, \|v^{k-1}P_-(v_{x})\|_{L^1_tL^2_x} \nonumber \\
               & \lesssim & \| v\|_{L^{\infty}_tH^1_x} \|v^{k-1}\|_{L^\infty_{t,x}} \|v_x\|_{L^1_tL^2_x}\nonumber \\
               & \lesssim & T \, \| v\|_{L^{\infty}_tH^1_x}^{k+1} \nonumber \\
               & \lesssim & T \, \| v\|_{X^1_{T,\lambda}}^{k+1} \, .
\end{eqnarray}
Next from obvious calculation, the Sobolev embeding (in space) $H^1
\hookrightarrow L^\infty$  and  Lemma \ref{lemderiv} we infer that
\begin{eqnarray}
   \|\partial_x c\|_{L^1_tL^2_x} & \lesssim & \| M(v^k) \, v \, M[v^{k-1} P_-(v_{x})]\|_{L^1_tL^2_x}
                                  + \| v_x M[v^{k-1} P_-(v_{x})]\|_{L^1_tL^2_x}\nonumber \\
                                  &+&  \| v^{k-1} \, v_x P_-(v_x)\|_{L^1_tL^2_x}
                                  + T^{1/2} \, \| P_+ [ e^{-iF}  \,  v^{k} P_-(v_{xx}) ] \|_{L^2_tL^2_x}\nonumber \\
                                & \lesssim &   \| v \|_{X^1_{T,\lambda}}^{2k} \, \|v_x\|_{L^1_tL^2_x}
                                 + \| v\|_{X^1_{T,\lambda}}^{k-1} \| v_{x}\|_{L^2_tL^4_x}^2\nonumber \\
                                 &+&T^{1/2} \, \| \partial_x(e^{-iF} v^{k} ) \|_{L^4_{t,x}} \,  \| P_-(v_{x}) \|_{L^4_{t,x}}\nonumber \\
                                & \lesssim & T \, \| v \|_{X^1_{T,\lambda}}^{2k+1}
                                 + T^{1/2} \| v\|_{X^1_{T,\lambda}}^{k+1} \, (1+\| v\|_{X^1_{T,\lambda}}^{k}) \, .
\end{eqnarray}

\noindent {\bf .} Estimates for $d=-k(k-1) P_+[vhe^{-iF}]$ where
\begin{equation}
   h=\sum_{q \in {\Z^* / \lambda} }\frac{C_q(v^{k-2}v_xHv_x)}{iq} \, e^{iqx} \, .
\end{equation}
We note first that,
\begin{eqnarray}
   \|d\|_{L^1_tL^2_x} & \lesssim & \| v \|_{L^1_tL^2_x}\, \| h\|_{L^\infty_{t,x}} \nonumber \\
               & \lesssim &  T \, \| v \|_{L^{\infty}_tL^2_x} \, \|v^{k-2}v_xHv_x\|_{L^\infty_tL^1_x} \nonumber \\
               & \lesssim & T \, \| v\|_{X^1_{T,\lambda}} \,  \|v^{k-1}\|_{L^\infty_{t,x}} \|v_x\|_{L^\infty_tL^2_x}^2 \nonumber \\
               & \lesssim & T \, \| v\|_{X^1_{T,\lambda}}^{k+1} \, .\nonumber \\
\end{eqnarray}
On the other hand, in the same way than previously,  we see that
\begin{eqnarray}
   \|\partial_x d\|_{L^1_tL^2_x} & \lesssim & \| v_x h \|_{L^1_tL^2_x} + \| M(v^k) v h \|_{L^1_tL^2_x}
                                     + \| v \sum_{q \in {\Z^* / \lambda}} C_q(v^{k-2}v_xHv_x) e^{iqx}\|_{L^1_tL^2_x} \nonumber \\
                                 & \lesssim &   \| v_x \|_{L^1_tL^2_x} \,
                                 \|v^{k-2}v_xHv_x\|_{L^\infty_tL^1_x} \nonumber \\
                                 &+& \| M(v^k) \|_{L^\infty_{t,x}} \|v \|_{L^1_tL^2_x} \| v^{k-2}v_xHv_x\|_{L^\infty_tL^1_x}\nonumber \\
                                 &+& T^{1/2} \| v \sum_{q \in {\Z^* / \lambda}} C_q(v^{k-2}v_xHv_x) e^{iqx} \|_{L^2_{t,x}} \nonumber \\
                                 & \lesssim & T \, \| v\|_{X^1_{T,\lambda}}^{k+1} + T \, \| v\|_{X^1_{T,\lambda}}^{2k+1}
                                 + T^{1/2} \| v \|_{L^\infty_{t,x}} \, \| v^{k-2}\|_{L^\infty_{t,x}} \| v_xHv_x\|_{L^2_{t,x}} \nonumber \\
                                 &  \lesssim & T \, \| v\|_{X^1_{T,\lambda}}^{k+1} + T \, \| v\|_{X^1_{T,\lambda}}^{2k+1}
                                 + T^{1/2} \| v \|_{L^\infty_{t,x}}^{k-1}\,  \| v_x\|_{L^4_{t,x}}^2 \nonumber \\
                                 & \lesssim & (T +T^{1/2}) \, \| v\|_{X^1_{T,\lambda}}^{k+1} + T \, \|v\|_{X^1_{T,\lambda}}^{2k+1}\,  .\nonumber \\
\end{eqnarray}

\section{Local well-posedness for $(GBO)$}
\subsection{Local well-posedness for "small" initial data}
\subsubsection{Nonlinear Estimates on $ u $}
As for the $(BO)$ equation, we first state the result for "small"
data. More precisely, for any real number $ A>1 $ given, we first
prove the local well-posedness
 result for $ u_0\in H^1_\lambda $ such that
\begin{equation}
 \|u_{0}\|_{L^2_\lambda}  \lesssim A \, \varepsilon^{1/k-1/2}   \quad \mbox{ and } \quad
 \|\partial_x u_0\|_{L^2_\lambda}  \lesssim \varepsilon^{1/2+1/k}  \label{d1}
 \end{equation}
 for some $ 0<\varepsilon=\varepsilon(A)<\!\! 1 $ which does not depend on $ \lambda $. Note that the $ L^2 $-norm of $ u_0 $
 may be taken  arbitrary large in \re{d1}. We stress out the attention of the reader that this will be necessary to consider such initial data
 since the $ L^2 $-norm is  subcritical
 for $(GBO)$ as soon as $k\geq 2$ and since we will use some
 dilation arguments in the case of non small initial data.

 So let us consider $ u_0\in H^\infty_\lambda $ satisfying \re{d1} and let $ u $ be the  emanating maximal solution given for instance by
  \cite{ABFS}. We are going to prove that $ u\in C([0,1];H^\infty_\lambda) $ with
\begin{equation}
\|u\|_{X^0_{1,\lambda}} \lesssim  A \, \varepsilon^{1/k-1/2} \quad
\mbox{ and } \quad \|\partial_x u\|_{X^0_{1,\lambda}} \lesssim
\varepsilon^{1/2+1/k} \label{d2} \quad .
\end{equation}

Notice that despite the local existence is only known for
$H^2_{\lambda}$ initial data, it will be enough to derive
$X^1_{T,\lambda}$ estimates for $u$ since we will check at the end
of this subsection that for $ 0<T\le 1 $, $ \|u\|_{L^\infty_T
H^2_\lambda} $-norm cannot blow up  as long as  $
\|u\|_{X^1_{T,\lambda}} $ remains bounded. Therefore, according to
the local well-posedness result in \cite{ABFS}, the solution can be
extended  in $ C([0,T];H^\infty_T) $ as soon as
$\|u\|_{X^1_{T,\lambda}} <\infty $.

 First, by a continuity argument, we can assume that
 $ \|\partial_x u \|_{X^0_{T,\lambda}} \lesssim \varepsilon^{1/2} $ for some $ 0<T<T^* $ where $ T^* $ is the maximal time of existence of the solution $ u $. Next, since the $ L^2 $-norm of $ u $ is a constant of the motion, to prove the desired result, it suffices to prove that if
\begin{equation}
\|u\|_{L^\infty_{T} L^2_\lambda} \lesssim A \,
\varepsilon^{1/k-1/2} \quad \mbox{ and } \quad \|\partial_x u
\|_{X^0_{T,\lambda}} \lesssim \varepsilon^{1/2}, \label{d3}
\end{equation}
with $ 0<T< 1 $, then
\begin{equation}
  \|u\|_{L^4_{T,\lambda}} \lesssim A \,
\varepsilon^{1/k-1/2} \quad \mbox{ and } \quad
 \|\partial_x u \|_{X^0_{T,\lambda}}\lesssim \varepsilon^{1/2+1/k} \quad .\label{d4}
\end{equation}
The estimate on $ \|u\|_{L^4_{1,\lambda}} $ is trivially satisfied since by Sobolev inequalities and interpolation,
$$
\|u\|_{L^4_{1,\lambda}}\lesssim \|u\|_{L^\infty_1\dot{H}^{1/4}_\lambda}
\lesssim \|u\|_{L^\infty_1L^2_\lambda}^{3/4}  \, \|u_x\|_{L^\infty_1 L^2_\lambda}^{1/4}
 \lesssim \varepsilon^{3/(4k) - 1/4} \lesssim \varepsilon^{1/k - 1/2} \quad .
$$
 We now estimate  the $ L^\infty_T L^2_\lambda $-norm of $ u_x $. Recall that since $u$ is real valued,
\begin{equation}
   \| u_x\|_{L^\infty_T L^2_\lambda} \leq  \| P_1(u_x) \|_{L^\infty_T L^2_\lambda} +  2\, \| P_+(u_x)\|_{L^\infty_T L^2_\lambda}
   \label{ppbis}
\end{equation}
with
\begin{equation}
   P_+(u_x)=P_+(\partial_x(e^{iF} w) )+P_+\Bigl( \partial_x  [e^{iF} P_-(e^{-iF} u)]  \Bigr)  \; , \label{pp}
\end{equation}
We consider the first term in the right hand side of \re{ppbis}. Since $u$ solves $(GBO)$, it follows from Bernstein inequalities that,
\begin{eqnarray}
  \| P_1 u_x\|_{L^\infty_T L^{2}_\lambda}
   & \lesssim & \|V(t) \partial_x u_{0} \|_{L^2_\lambda}
   + \Bigl\| \int_0^t V(t-t') u^{k+1}(t') \, dt' \Bigr\|_{L^\infty_T L^{2}_\lambda} \nonumber \\
   & \lesssim & \| \partial_x u_{0} \|_{L^2_\lambda }
   + T \, \|u\|_{L^\infty_T L^{2(k+1)}_\lambda}^{k+1} \nonumber \\
   & \lesssim & \| \partial_x u_{0} \|_{L^2_\lambda }
   + T \, \|u\|_{X^1_{T, \lambda }}^{k+1} \nonumber \\
   & \lesssim & \varepsilon^{1/k+1/2} \,
   \label{d9bis}
\end{eqnarray}
choosing $T=T(\varepsilon)>0$ small enough.

\noindent We consider now the second  term in the right hand side of
\re{ppbis} by means of the decomposition given by \re{pp}. From
\re{estprop} in Proposition 4.1 we infer that for
$T=T(\varepsilon)>0 $ small enough,
\begin{eqnarray}
  \|\partial_x P_+ (e^{iF} w)\|_{L^\infty_T L^2_\lambda}  & \lesssim & \| F_x\,  w \|_{L^\infty_T L^2_\lambda}
 + \|\partial_x w \|_{L^\infty_T L^2_\lambda} \nonumber \\
  & \lesssim & \|u\|^{k+1}_{L^\infty_T L^{2(k+1)}_\lambda} +\|\partial_x w_0\|_{L^2_\lambda}+
  T^{1/4} \, \|u\|^{3k+1}_{X^1_{T,\lambda}} \nonumber \\
  & \lesssim &  \|u\|^{k+1}_{L^\infty_T L^{2(k+1)}_\lambda}
  +\|u_{0}\|^{k+1}_{L^{2(k+1)}_\lambda}   +  \|\partial_x u_{0}\|_{L^2_\lambda}+  \varepsilon ^{1/2+1/k}\nonumber \\
  &\lesssim &  \|u\|^{k+1}_{L^\infty_T L^{2(k+1)}_\lambda}
  +\|u_{0}\|^{k+1}_{L^{2(k+1)}_\lambda}
  +  \varepsilon ^{1/2+1/k}\label{d5} \; .
\end{eqnarray}
Moreover, using Lemma 4.1 we see that
\begin{eqnarray}
 \Bigl\|\partial_x P_+\Bigl( e^{iF} \, P_-(e^{-iF}  u)\Bigr) \Bigr\|_{L^\infty_T L^2_\lambda}
 & \lesssim &\|\partial_x F \|_{L^\infty_T L^{2(k+1)/k}_\lambda}  \, \|u \|_{L^\infty_T L^{2(k+1)}_\lambda} \nonumber \\
 & \lesssim &  \|u \|_{L^\infty_T L^{2(k+1)}_\lambda}^{k+1}
 \label{d6}\quad .
\end{eqnarray}
It thus remains to get a good estimate on  $ \|u_0 \|_{ L^{2(k+1)}_\lambda}^{k+1} $ and  $ \|u \|_{L^\infty _T L^{2(k+1)}_\lambda}^{k+1} $. To do this we first note that we have,
\begin{eqnarray}
   \|u_0 \|_{ L^{2(k+1)}_\lambda}^{k+1}
  &\lesssim & \| u_0 \|_{\dot{H}^{k/[2(k+1)]}_\lambda}^{k+1} \nonumber \\
  &\lesssim & \| u_0\|_{\dot{H}^1_\lambda}^{k/2} \, \| u_0\|_{L^2_\lambda}^{(k+2)/2} \nonumber \\
  &\lesssim & \varepsilon ^{(1/2+1/k)k/2} \, \varepsilon ^{(1/k-1/2)(k+2)/2} \nonumber \\
  &\lesssim & \varepsilon ^{1/2+1/k}  \;  .
\label{d8}
\end{eqnarray}
Concerning $ \|u \|_{L^\infty_T L^{2(k+1)}_\lambda}^{k+1} $, remark that,
\begin{eqnarray}
   \|u \|_{L^\infty_T L^{2(k+1)}_\lambda}^{k+1}
   &\lesssim &\|P_1 u\|_{L^\infty_T L^{2(k+1)}_\lambda}^{k+1}+
   \|\partial_x P_{>1} u \|_{L^\infty_TL^2_\lambda}^{k+1} \nonumber \\
  &\lesssim & \|P_1 u\|_{L^\infty_T L^{2(k+1)}_\lambda}^{k+1}+ \varepsilon ^{(k+1)/2} \; .
\label{d8}
\end{eqnarray}
But using that $u$ is a solution of $(GBO)$, Sobolev, Bernstein  and Strichartz estimates, \re{d1} and  \re{d3}, we infer that
\begin{eqnarray}
  \|P_1 u \|_{L^\infty_T L^{2(k+1)}_\lambda}^{k+1}
   & \lesssim & \|V(t) u_{0} \|_{\dot{H}^{k/[2(k+1)]}_\lambda}^{k+1}
   + \Bigl\| \int_0^t V(t-t') u^{k+1}(t') \, dt' \Bigr\|_{L^\infty_T L^{2}_\lambda}^{k+1} \nonumber \\
   & \lesssim & \|u_{0} \|_{\dot{H}^{k/[2(k+1)]}_\lambda}^{k+1}
   + T^{k+1} \, \|u\|_{L^\infty_T L^{2(k+1)}_\lambda}^{(k+1)} \nonumber \\
   & \lesssim & \varepsilon^{1/k+1/2} + T^{k+1} \, \|u\|_{L^\infty_T L^{2(k+1)}_\lambda}^{(k+1)} \,
   \label{d9}
\end{eqnarray}
since the operators $\partial_x \, P_1$ and $V(\cdot)$ are respectively  bounded from $L^{2(k+1)}$  to $L^2$ and from $L^2$ to itself.
This ensures that for $T=T(\varepsilon )>0$ small enough,
\begin{equation}
  \|P_1 u\|_{L^\infty_T L^{2(k+1)}_\lambda}^{k+1} \lesssim
  \varepsilon^{1/k+1/2} \quad .
  \label{d9bis}
\end{equation}
Gathering \re{ppbis}-\re{d9bis}, we thus obtain for $ 0<T<1 $,
\begin{eqnarray}
  \|\partial_x u \|_{L^\infty_T L^2_\lambda} & \lesssim & \varepsilon^{1/2+1/k} \; .
   \label{d14}
\end{eqnarray}
Let us now  estimate the $ L^4_{T,\lambda} $-norm of $  u_x $.
Proceeding in the same way than previously we first infer that
\begin{eqnarray}
  \| P_1 u_x\|_{L^4_{T,\lambda}}
   & \lesssim & \|V(t) \partial_x u_{0} \|_{L^4_{T,\lambda}}
   + \Bigl\| \int_0^t V(t-t') u^{k+1}(t') \, dt' \Bigr\|_{L^4_{T,\lambda }} \nonumber \\
   & \lesssim & \| \partial_x u_{0} \|_{L^2_\lambda }
   + T \, \|u\|_{L^\infty_T L^{4(k+1)}_\lambda} \nonumber \\
   & \lesssim & \varepsilon^{1/k+1/2} \,
   \label{d9bis}
\end{eqnarray}
choosing $T=T(\varepsilon)>0$ small enough.

\noindent Next we have,
  \begin{eqnarray}
 \|\partial_xP_+ (e^{iF} w)\|_{L^4_{T,\lambda}}  & \lesssim & \| F_x\,  w \|_{ L^4_{T,\lambda}}
 + \|\partial_x w \|_{ L^4_{T,\lambda}} \nonumber \\
  & \lesssim &  \|u\|^{k+1}_{ L^{4(k+1)}_{T,\lambda}}
  +\|u_{0}\|^{k+1}_{L^{2(k+1)}_\lambda}  \nonumber \\
  &+&  \|\partial_x u_{0}\|_{L^2_\lambda}+
 T^{1/4} \|u \|^{k+1}_{X^1_{T,\lambda}} \label{d15} \quad .
 \end{eqnarray}

\noindent This combines with \re{d14} and \re{pp} completes the proof of
\re{d2}, since by \re{d1}, \re{d14} and Sobolev inequalities,
\begin{eqnarray}
     \|u \|_{L^{4(k+1)}_{T,\lambda}}^{k+1} &\lesssim & \|u \|_{L^\infty_T\dot{H}^{(2k+1)/[4(k+1)]}_\lambda }^{k+1} \nonumber \\
      &\lesssim & \|u \|_{L^\infty_T\dot{H}^{1}_\lambda}^{(2k+1)/4} \, \|u \|_{L^\infty_TL^2_\lambda }^{(2k+3)/4} \nonumber \\
      & \lesssim & \varepsilon ^{(1/2+1/k)(2k+1)/4}  \, \varepsilon ^{(1/k-1/2)(2k+3)/4} \nonumber \\
      &  \lesssim & \varepsilon ^{3/4+1/k} \; .
\end{eqnarray}

It remains to check that $ \|u_\lambda\|_{X^2_{T,\lambda}} $ can
not go to infinity as long as  $\|u_\lambda\|_{X^1_{T,\lambda}} $
remains bounded. First we have (here we use \re{estbisprop}),
\begin{eqnarray}
     \|\partial^2_x P_+ (e^{-iF} w)\|_{X^0_{T,\lambda}}  & \lesssim &
     \| \partial_x F \,  \partial_x w\|_{X ^0_{T,\lambda}}+\|\partial^2_x w\|_{X^0_{T,\lambda}}\nonumber \\
     &+& \| \partial_x^2 F \, w \|_{X^0_{T,\lambda}} +\| (\partial_x F)^2 \,  w\|_{X^0_{T,\lambda}}\nonumber \\
     &  \lesssim & \|\partial^2_x w \|_{X^0_{T,\lambda}}+ \|u\|_{L^\infty_{T,\lambda}}^k   (1+\|u\|_{L^\infty_{T,\lambda}}^k)    \| w \|_{X^1_{T}} \nonumber \\
      &+& \|u\|_{L^\infty_{T,\lambda}}^{k-1} \|w\|_{L^\infty_{T,\lambda}} \, \|\partial_x u \|_{X^0_{T,\lambda}} \nonumber \\
     & \lesssim & \|u\|_{X^1_{T,\lambda}}^{k}(1+\|u_\lambda\|_{X^1_{T,\lambda}}^k) \|w\|_{X^1_{T,\lambda}} \nonumber \\
     &+& \|\partial_x^2 w_{0} \|_{L^2_\lambda} +T^{1/4} \,  Q(\|u\|_{X^1_{T,\lambda}}) \|u\|_{X^2_{T,\lambda}} \,. \label{d20bis}
\end{eqnarray}
Furthermore,
\begin{eqnarray}
   \Bigl\|\partial_x^2  P_+ \Bigl( e^{iF} \, P_-(e^{-iF} u)\Bigr)  \Bigr\|_{X^0_{T,\lambda}}
   & \lesssim &   (\|\partial_x^2 F\|_{L^\infty_{T} L^2_\lambda}   +\|\partial_x F\|_{L^\infty_{T} L^4_\lambda}^2)\|u\|_{L^\infty_{T,\lambda}}
    \nonumber \\
   & \lesssim &  \|u\|_{X^1_{T,\lambda}}^{k+1}(1+ \|u\|_{X^1_{T,\lambda}}^{k}) \quad . \label{d21bis}
\end{eqnarray}
Hence, gathering \re{d20bis}-\re{d21bis}, we infer that
\begin{equation}
   \|u\|_{X^2_{T,\lambda}} \lesssim   \|u_{0}\|_{H^2}(1+\|u_{0}\|_{H^2}^{2k})
   +\|u\|_{X^1_{T,\lambda}}(1+\|u\|_{X^1_{T,\lambda}}^{2k+1})
   +T^{1/4}  Q(\|u\|_{X^1_{T,\lambda}}) \, \|u\|_{X^2_{T,\lambda}} \, .\nonumber
\end{equation}
This completes the proof of \re{d2}.
\subsubsection{Local existence}
Now, let $ u_0\in H^{1}_{\lambda} $ satisfying \re{d1} and let $\{
u_{0}^n \} \subset H^{\infty}_{\lambda} $ a sequence converging to $
u_0 $ in
 $ H^{1}_{\lambda} $. We denote by $ u_n $ the solution of $(GBO) $   emanating from $ u_{0}^n $.
 From the previous results  $ u_n\in C([0,1];H^\infty_\lambda) $ and moreover,
   $ \|u_n \|_{X^1_{1,\lambda}} \lesssim \varepsilon^{1/k-1/2} $ uniformly in $ n $. Thus we can pass to the limit up to a subsequence which leads to the existence of a solution $ u\in X^{1}_{1,\lambda} $ of $(GBO)$ with $ u_0 $ as initial data.
\subsubsection{ Continuity, uniqueness and  regularity of the flow map }
As for the $(BO)$ equation  we   first prove  that the flow-map is Lipschitz on a
small ball of $ H^1_\lambda $. The continuity of $ t\mapsto u (t)
$ in $ H^1_\lambda $ will follow directly.

Let  $ u_1 $ and $u_2 $ be two solutions of $(GBO)$ in $ X^1_{T,\lambda} $,  associated
with the initial data $ \varphi_1 $ and $ \varphi_2 $ in $ H^1_\lambda $. We assume that they satisfy
\begin{equation}
\|u_i\|_{X^1_{T,\lambda}} \lesssim A \, \varepsilon^{1/k-1/2}
\quad \mbox{ and } \quad
\|\partial_x u_i\|_{X^0_{T,\lambda}} \lesssim \varepsilon^{1/2+1/k} \; , i=1,2 \label{d20} \quad ,
\end{equation}
where $ 0<\varepsilon=\varepsilon(A)<\!\! 1 $ has the same value  as in \re{d1}.
 We then consider
 $$
 z=w_1-w_2=-i P_+(e^{-iF_1} u_1)+iP_+ (e^{-iF_2}
  u_2)
 $$
 with $ F_i $ defined as $ F $ with $ u $ replaced by $ u_i $.
 Following the calculations performed in Subsection 4.2 we clearly have,
 \begin{equation}
 \|z\|_{X^1_{T,\lambda}}  \lesssim  \| \varphi_1-\varphi_2\|_{H^1}(1+\|\varphi_1\|_{H^1}^k)
  +T^\nu \, \|z\|_{X^1_{T,\lambda}} (\|u_1\|_{X^1_{T,\lambda}}^k +\|u_2\|_{X^1_{T,\lambda}}^k) \label{d19}
 \end{equation}
 Setting
 $$
\Lambda=\sum_{i=1}^2 \|u_i\|_{L^\infty_T L^{2(k+1)}_\lambda}
 + \|u_i\|_{L^{4(k+1)}_{T,\lambda}} \quad ,
 $$
 clearly, $ \|w_i\|_{L^\infty_T L^{2(k+1)}_\lambda}
 + \|w_i\|_{L^{4(k+1)}_{T,\lambda}}\lesssim  \Lambda $, and  by
 \re{d2} and Sobolev inequalities,
 \begin{equation}
 \Lambda \lesssim \varepsilon^{1/k-\frac{1}{2(k+1)}}\label{yoyo} \quad .
 \end{equation}
 From \re{pp} and Lemmas 4.1, we get after straightforward computations,
\begin{eqnarray}
\|\partial_x (u_1 -u_2) \|_{X^0_{T,\lambda}} & \lesssim & \|
\partial_x z\|_{X^0_\lambda} (1+\Lambda^k) +\Lambda^{k}\| u_1 -u_2\|_{X^1_{T,\lambda}}\nonumber \\
  & & + \|(e^{iF_1}-e^{iF_2})\|_{L^\infty_{T,\lambda}}
 \Bigl( \|\partial_x w_2\|_{X^0_{T,\lambda}} +
\Lambda^{k+1}\Bigr) \; . \label{d21}
\end{eqnarray}
 On the other hand,
$$
    \|F_1-F_2|\|_{L^\infty_{T,\lambda }}  \lesssim \|u_1^k-u_2^k \|_{L^\infty_T L^1_\lambda} \lesssim \|u_1-u_2\|_{L^\infty_T L^k_\lambda}
    (\|u_1\|_{L^\infty_T L^k_\lambda}+\|u_2\|_{L^\infty_T L^k_\lambda})^{k-1} \; .
$$
Hence, noticing that by Sobolev inequality,
$$
\|u_i\|_{L^\infty_T L^k_\lambda}\lesssim \|u_i\|_{L^\infty_T
L^2_\lambda}^{1/2+1/k} \|\partial_x u_i\|_{L^\infty_T
L^2_\lambda}^{1/2-1/k} \lesssim 1 , \quad i=1,2 \quad ,
$$
we thus obtain that,
\begin{equation}
   \|e^{iF_1}-e^{iF_2}\|_{L^\infty_{T,\lambda}}
   \lesssim  \| F_1-F_2 \|_{L^\infty_{T,\lambda}}
   \lesssim \|u_1-u_2\|_{X^1_{T,\lambda}}  \; .
\end{equation}
Noticing also that
$$
    \|\partial_x w_2 \|_{X^0_{T,\lambda}}
    \lesssim \|\partial_x u_2\|_{X^0_{T,\lambda}} +\Lambda^{k+1}
   \lesssim \varepsilon^{1/2+1/k}  \; ,
 $$
and
\begin{equation}
\|(u_1-u_2)\|_{X^0_{T,\lambda}} \lesssim \|u_1-u_2\|_{L^2_\lambda}+ T^{1/4} \|u_1-u_2\|_{X^1_{T,\lambda}}
(\|u_1\|_{X^1_{T,\lambda}}+\|u_2\|_{X^1_{T,\lambda}}) \label{d22}
\end{equation}
we finally deduce from \re{d19}-\re{d22} that
\begin{equation}
   \|u_1-u_2 \|_{X^1_{T,\lambda}}
   \lesssim   \|\varphi_1-\varphi_2\|_{H^1_\lambda} (1+\|\varphi_1\|_{H^1_\lambda}) \quad . \label{d23}
\end{equation}
 Combining \re{d23} with the same arguments as in the end of Section 3.2, we obtain that the solution $ u$ constructed above
 belongs to $ X^1_{1,\lambda}\cap C([0,1],H^1_\lambda) $ and is unique in this class. Moreover,  the flow-map is Lipschitz from
  the ball of $ H^1_\lambda $  $$ \{\varphi\in H^1_\lambda, \quad \|\varphi\|_{L^2_\lambda} \le A \varepsilon^{1/k-1/2}, \,
   \|\varphi_x \|_{L^2_\lambda} \le \varepsilon^{1/2+1/k} \, \} ,  \, \mbox{ with }\varepsilon=\varepsilon(A),$$
    into $X^1_{1,\lambda} \cap C([0,1],H^1_\lambda) $

\subsection{Arbitrary large initial data}
We used the dilation symmetry argument to extend the result for
arbitrary large data. First note that if $ u(t,x) $ is a
$2\pi$-periodic solution of $(GBO)$ on $[0,T] $ with initial data $
u_0 $, then $ u_\lambda(t,x)=\lambda^{-1/k}
u(\lambda^{-2}t,\lambda^{-1} x)$ is a $ 2\pi\lambda $-periodic
solution of $(GBO)$ on
 $ [0,\lambda^2 T] $ emanating from the initial data $u_{0,\lambda}=\lambda^{-1/k} u_0(\lambda^{-1} x) $.

Let $ u_0\in H^1 $.  If $ \|\partial_x u_0\|_{L^2} \le \varepsilon^{1/2+1/k} $ with
$$
   \varepsilon=
  \varepsilon [(\|\partial_x u_0\|_{L^2}^\frac{k-2}{k+2}+1)\|u_0\|_{L^2}]<1 \; ,
$$
then $ u_0 $ satisfies \re{d2} with $ A=(\|\partial_x
u_0\|_{L^2}^\frac{k-2}{k+2}+1)\|u_0\|_{L^2} $ and so we are done.  Otherwise, we set
$$
\lambda=\varepsilon^{-1} \|\partial_x u_0\|_{L^2}^\frac{2k}{k+2} \ge 1 \; ,
$$
 so that $ u_{0,
\lambda} $ satisfies \re{d2} with $ \varepsilon $ and $ A $ defined
as above. We are thus reduced to the case of small initial data.
Therefore, there exists a unique local solution $ u_\lambda\in
C([0,1] , H^1_\lambda)\cap X^1_{1,\lambda} $ of $(GBO)$  emanating
from $ u_{0,\lambda} $. This proves the existence
 and uniqueness in $ C([0,T] , H^1)\cap X^1_T $ of the solution $ u $ emanating from $ u_0 $ with $ T=T(\|u_0\|_{H^1})
  $ and $ T\rightarrow +\infty $ as $ \|u_0\|_{H^1} \rightarrow 0 $. The fact that the flow-map is Lipschitz on every bounded set of $ H^1 $ follows as well.

 Finally, note that the change of unknown \re{chgtvarGBO} preserves the continuity of the solution and the continuity
 of the flow-map in $ H^1 $. Moreover, for $k=2$, the flow-map is Lipschitz  on every  closed set $ S_\beta $ of $ H^1(\T) $ of the form
 $$
 S_\beta=\{\varphi \in H^1(\T) , \quad \into \varphi^2 =\beta \quad \} \quad .
 $$
 On the other hand, it does not preserve the Lipschitz property of the flow. Therefore,
  contrary to   the real line case (cf.\cite{MR2}), we do not know if the flow-map is Lipschitz or even uniformly continuous
  on bounded set. Recall that on the real line, the flow-map is known to be real-analytic on a small ball
   of $ H^1(\R) $ (cf. \cite{KPV} and \cite{MR1}) and Lipschitz on every bounded set of $ H^1(\R) $ (cf. \cite{MR2}).


\begin{thebibliography}{10}
\bibitem{ABFS} {\sc L. Abdelouhab, J. Bona, M. Felland, and J.C. Saut},
Nonlocal models for nonlinear, dispersive waves,
{\em Phys. D {\bf 40} (1989),  360--392.}

\bibitem{AF} {\sc M.J. Ablowitz, A.S. Fokas},
 The inverse scattering transform for the Benjamin- Ono equation, a pivot for multidimensional problems,
{\em Stud. Appl. Math. 68 (1983), 1--10.}

\bibitem{B} {\sc T.B. Benjamin},
Internal waves of permanent form in fluids of great depth,
{\em J. Fluid Mech. {\bf 29} (1967), 559--592.}

\bibitem {BL} {\sc H.A.~Biagioni  and F.~Linares},
Ill-posedness for the derivative Schr\"odinger and generalized Benjamin-Ono equations,
{\em  Trans. Amer. Math. Soc. {\bf 353} (2001),  3649-3659.}

\bibitem {BoKa}  {\sc J. Bona and H. Kalisch},  Singularity formation in the generalized Benjamin-Ono equation,
{\em To appear in Disc. Cont. Dyn. Systems-Ser. B. }

\bibitem {Bo1}{\sc J.~Bourgain},
Fourier transform restriction phenomena for certain lattice subsets and application to nonlinear evolution equations II. The Schr\"odinger equation,
{\em GAFA, {\bf 3}  (1993), 209-262}

\bibitem{CKSTT} {\sc J.~Colliander, M.~Keel, G.~Staffilani, H.~Takaoka, T.~Tao},
Sharp global well-posedness results for periodic and non-periodic KdV and modified KdV on $\R$ and $\T$,
{\em JAMS {\bf16}  (2003),   705-749.}

\bibitem{CW} {\sc R.~Coifman, M.~Wickerhauser},
{\it The scattering transform for the Benjamin-Ono equation},
Inverse Probl. 6 (1990), 825-860.

\bibitem{GV2} {\sc J. Ginibre and G. Velo},
Smoothing properties and existence of solutions for the generalized Benjamin-Ono equation,
{\em J. Differential Equations {\bf 93} (1991), 150--232.}

\bibitem {HaOz} {\sc N.~Hayashi and T.~Ozawa},
Remarks on Schr\"odinger equations in one space dimension,
{\em Diff. Int. Equ. {\bf 7} (2) (1994),  453-461.}

\bibitem {Io} {\sc J.R.~Iorio},
On the Cauchy problem for the Benjamin-Ono equation,
{\em Comm. Partial Differential Equations {\bf 11} (10) (1986),  1031-1081.}

\bibitem {IK} {\sc A.D.~Ionescu, C.E.~Kenig},
Global well posedness of the Benjamin-Ono equation in low regularity
spaces, {\em Preprint (2005).}

\bibitem {KK} {\sc C.E.~Kenig and K.~Koenig},
On the local well-posedness of the  Benjamin-Ono and modified Benjamin-Ono equations,
{\em Preprint (2005).}

\bibitem {KPV2} {\sc C.E.~Kenig, G.~Ponce and L.~Vega},
Well-posedness of the initial value problem for the Korteweg-de Vries equation,
{\em J. Amer. Math. Soc. {\bf 4} (2) (1991),  323-347.}

\bibitem {KPV} {\sc C.E.~Kenig, G.~Ponce and L.~Vega},
On the Generalized Benjamin-Ono equations,
{\em  Trans. Amer. Math. Soc. {\bf 342} (1994), 155-172.}

\bibitem {KT1} {\sc H.~Koch and N.~Tzvetkov},
On the local well-posedness of the Benamin-Ono equation in $H^s(\R)$,
{\em IMRN  {\bf 26} (2003), 1449-1464.}

\bibitem {KT2} {\sc H.~Koch and N.~Tzvetkov},
Nonlinear wave interactions for  the Benamin-Ono equation,
{\em Preprint (2002).}


\bibitem{MR1}{\sc L.~Molinet and F.~Ribaud},
 Well-posedness results for the Generalized Benjamin-Ono equation with small initial data,
{\em J. Math. Pures Appl. {\bf 83} (2004), 277-311.}

\bibitem{MR2}{\sc L.~Molinet and F.~Ribaud},
 Well-posedness results for the generalized Benjamin-Ono equation with arbitrary large initial data,
 {\em I.M.R.N. {\bf 70}, (2004), p.p. 3757-3795.}

\bibitem {MST} {\sc L.~Molinet, J.C.~Saut and N.~Tzvetkov},
Ill-posedness issues for the Benjamin-Ono and related equations,
{\em SIAM J. Math. Anal. {\bf 33} (4) (2001),  982-988.}

\bibitem {Po}{\sc G.~Ponce},
On the global well-posedness of the Benjamin-Ono equation,
{\em Differential Integral Equations  {\bf 4} (3) (1991),  527-542.}

\bibitem {Saut}{\sc J.C.~Saut},
Sur quelques generalisation de l' equation de Korteweg-de Vries,
{\em J. Math. Pures Appl. {\bf 58} (1979),  21-61.}

\bibitem {Tao} {\sc T.~Tao},
Global well-posedness of the Benjamin-Ono equation in $ H^1(\R) $,
{\em Preprint (2003).}
\end{thebibliography}
\end{document}